\documentclass[11pt,reqno]{amsart}

\usepackage{euscript}
\usepackage{epsfig}
\usepackage{amssymb}
\usepackage{epic}

\theoremstyle{plain}

\theoremstyle{definition}

\newcommand{\bms}{\mbox{\boldmath$s$}}
\DeclareMathOperator{\Hom}{{\rm Hom}}
\DeclareMathOperator{\Ker}{{\rm ker}}
\DeclareMathOperator{\coker}{{\rm coker}}
\DeclareMathOperator{\rank}{{\rm rank}}
\newcommand{\ssw}{{\bf sw}}

\newcommand{\bH}{{\mathbb H}}

\def\C{\mathbb C}
\def\Q{\mathbb Q}
\def\R{\mathbb R}
\def\Z{\mathbb Z}
\def\N{\mathbb N}

\newcommand{\cale}{{\mathcal E}}

\newcommand{\calv}{{\mathcal V}}

\newcommand{\calt}{{\mathcal T}}

\newcommand{\calj}{{\mathcal J}}

\begin{document}

\title{On the Heegaard Floer homology of $S^3_{-p/q}(K)$}

\author{Andr\'as N\'emethi}
\address{Department of Mathematics\\Ohio State University\\Columbus, OH 43210, USA; and}
\address{R\'enyi Institute of Mathematics\\ Budapest, Hungary}
\address{}
\email{nemethi@math.ohio-state.edu; and nemethi@renyi.hu}
\urladdr{http://www.math.ohio-state.edu/\textasciitilde nemethi}
\thanks{The author is partially supported by NSF grant DMS-0304759.}



\begin{abstract}
Assume that the oriented 3-manifold  $M=S^3_{-p/q}(K)$ is
obtained by a rational  surgery (with coefficient $-p/q<0$) along
an algebraic knot $K\subset S^3$. We compute the Heegaard Floer
homology of $-M$ in terms of $p/q$ and the Alexander polynomial of
$K$.
\end{abstract}

\maketitle \pagestyle{myheadings} \markboth{{\normalsize Andr\'as
N\'emethi}}{ {\normalsize On the Heegaard Floer homology  of
$S^3_{-p/q}(K)$ }}

{\small

\section{Introduction}

In this article we compute the Heegaard Floer homology $HF^+(-M)$
(introduced by Ozsv\'ath and Szab\'o \cite{OSz}) for the oriented
3-manifold $M=S^3_{-p/q}(K)$ obtained by a negative rational
surgery (with coefficient $-p/q$) along an algebraic knot
$K\subset S^3$. In this case, since $H_1(M,\Z)=\Z_p$, the
$spin^c$-structures $\{\sigma_a\}_a$ of $M$ can be parametrized by
integers $a=0,1,\ldots p-1$. The main result of the article
establishes $HF^+(-M,\sigma_a)$ in terms of the integers $p$,
$q$, $a$, and the Alexander polynomial $\Delta$ of $K\subset S^3$.
Notice that the Alexander polynomial of an algebraic (any) knot
is well-understood, it can be easily computed from most of the
other invariants of the knot (e.g., in the present algebraic case,
from Puiseux or Newton pairs,
or from the semigroup associated with the corresponding local
analytic germ). In particular, the input of the theorem is the
simplest what one can hope. Since (in some sense) all
the coefficients of $\Delta$ are effectively involved in the description of
the Heegaard Floer homology, in fact, the result is optimal.

In the very  recent manuscript \cite{OSZuj}, Ozsv\'ath and Szab\'o
computed $HF^+(S^3_p(K))$ -- for any knot $K$ and any integer
surgery coefficient $p$ -- in terms of the filtered chain
homotopy type of the Heegaard Floer complex associated with the
pair $(S^3,K)$. Compared with this, our starting data,  and also
the description of $HF^+(-M)$, are simpler, and totally
elementary; as a price for this we have to impose the `negativity
restrictions' for the surgery coefficient and for $K$.

The proof (and the structure of the article) is based on the
results and constructions of  \cite{OSZINV} valid for plumbed
3-manifolds associated with some negative definite plumbing graphs
-- in fact, this also explains  the source of our restrictions.
Although \cite{OSZINV} presents a precise algorithm how one
should compute $HF^+$, its implementations in different
situations sometimes is not straightforward. In the present case
too, the proof and additional constructions run over many
sections. In fact, with the present article, we also wish to
advertise the efficiency, novelty and the power of \cite{OSZINV}.
This method, in fact, determines the `graded roots' (some graded
trees) associated with the plumbing graph of $M$, from  which one
can read  easily the Heegaard Floer homology. The advantage of
these graded roots is that (in all the cases known by the author)
their structure reflects perfectly the corresponding geometrical
construction which provides $M$. In particular, in many cases,
from the topology of $M$ one can  identify the roots rather
conceptually.

Section 2 recalls the classical invariants of algebraic knots and
connects them with the plumbing of $M$. The next section recalls
the definition and first properties of graded roots -- necessary
to formulate the main theorem, which appears in section 4.
Section 5 presents  two relevant results of \cite{OSZINV}, as
general principles to compute $HF^+$. (In fact, section 5 can
serve as a general recipe for the interested readers to compute
the Heegaard Floer homology in different other cases as well.)
The first result (which provides the absolute grading) is worked
out for our present situation in section 6, the other (regarding
the graded roots) in section 7. The last section contains some
examples.

\section{The manifold $M=S^3_{-p/q}(K_f)$}

\subsection{Review of algebraic knots}\label{alg1} \cite{BK,EN,GZDC}
Let $K_f\subset S^3$ be the link of an {\em irreducible complex
plane curve singularity } $f:(\C^2,0)\to (\C,0)$; i.e. for
$\epsilon>0$ sufficiently small, write $S^3=\{z\in \C^2:
||z||=\epsilon\}$ and take the transversal intersection
$K_f:=\{f=0\}\cap S^3\hookrightarrow S^3$. The natural
orientations of $S^3$ and of the regular part of $\{f=0\}$
induces a natural orientation on $K_f$.  Since $f$ is
irreducible, $K_f\approx S^1$. We will assume that $\{f=0\}$ is
not smooth at the origin,  i.e. $K_f\subset S^3$ is not the
unknot.  The isotopy type of $K_f\subset S^3$ is completely
characterized by any of the following invariants listed below.

\subsubsection{}\label{np} The {\em Newton pairs} of $f$ consist of $g\geq 1$ pairs
of integers $\{(p_i,q_i)\}_{i=1}^g$, where $p_i\geq 2$, $q_i\geq
1$, $q_1>p_1$ and gcd$(p_i,q_i)=1$.

\subsubsection{}\label{lp}   In some topological
constructions, it is preferable to replace the Newton pairs  by
the {\em `linking pairs'} (or, the decorations of the splice
diagram, cf. \cite{EN}) $(p_i,a_i)_{i=1}^g$, where
$$a_1=q_1 \ \mbox{ and} \  a_{i+1}=q_{i+1}+p_{i+1}p_ia_i \ \mbox{ for }\
i\geq 1.$$

\subsubsection{}\label{ap} Let $\Delta(t)$ be the {\em Alexander polynomial}
of $K_f\subset S^3$, or equivalently, the {\em characteristic
polynomial} of the algebraic monodromy acting on the first
homology of the Milnor fiber of $f$. It is normalized by
$\Delta(1)=1$. In terms of the $(p_i,a_i)_i$ pairs it is
$$\Delta(t)=\frac{ (t^{a_1p_1p_2\cdots p_g}-1)(t^{a_2p_2\cdots
p_g}-1)\cdots (t^{a_gp_g}-1)(t-1)}{(t^{a_1p_2\cdots
p_g}-1)(t^{a_2p_3\cdots p_g}-1)\cdots (t^{a_g}-1)(t^{p_1\cdots
p_g}-1)}.$$ The degree of $\Delta$, or equivalently, the first
Betti number of the Milnor fiber of $f$, is the Milnor number
$\mu$ of $f$.

\subsubsection{}\label{chp} The {\em semigroup $\Gamma $} of the germ
$f$ is a sub-semigroup of $\N$ defined by $\Gamma:= \{i_0(f,h):
h\in {\mathcal O}_{(\C^2,0)}\}$, where $i_0(f,h)$ denotes the
local intersection multiplicity of $f$ and $h$ (which equals the
codimension of the ideal $(f,h)$ in ${\mathcal O}_{(\C^2,0)}$).
For $h$ invertible $i_0(f,h)$ is zero, hence $0\in \Gamma$.

It is known that $\Gamma$ is generated by the integers
$\bar{\beta}_0=p_1p_2\cdots p_g$, $\bar{\beta}_k= a_kp_{k+1}\cdots
p_g$ for $1\leq k\leq g-1$, and $\bar{\beta}_g=a_g$. Moreover,
$\#(\N\setminus \Gamma)$ is finite. Its cardinality is the {\em
delta-invariant} $\delta$  of $f$, which in this case also equals
$\mu/2$, cf. \cite{Milnorbook}. It is also known that
$\delta=\mu/2$ equals the minimal Seifert genus of $K_f\subset
S^3$.

 The largest element of
$\N\setminus \Gamma$ is $\mu-1$.  In fact, for $0\leq k\leq
\mu-1$ one has: $k\in \Gamma$ if and only if $\mu-1-k\not\in
\Gamma$.

\subsubsection{}\label{resgr} The {\em embedded minimal good dual
resolution graph of the germ } $f$ has the following shape

\begin{picture}(400,80)(0,0)
\put(50,60){\circle*{4}} \put(100,60){\circle*{4}}
\put(150,60){\circle*{4}} \put(250,60){\circle*{4}}
\put(300,60){\circle*{4}} \put(100,20){\circle*{4}}
\put(150,20){\circle*{4}} \put(250,20){\circle*{4}}
\put(300,20){\circle*{4}} \put(370,60){\makebox(0,0){$K_f$}}
\dashline{3}(50,60)(175,60) \dashline{3}(225,60)(300,60)
\put(100,20){\dashbox{3}(0,40){}}
\put(150,20){\dashbox{3}(0,40){}}
\put(250,20){\dashbox{3}(0,40){}}
\put(300,20){\dashbox{3}(0,40){}}
\put(200,60){\makebox(0,0){$\cdots$}}
\put(300,60){\vector(1,0){30}} \put(300,70){\makebox(0,0){$-1$}}
\put(310,50){\makebox(0,0){$v_0$}}
\put(10,40){\makebox(0,0){$G(f):$}}
\end{picture}

\noindent
Above we emphasize only the vertices of degree one and three. The
dash-line between two such vertices replaces a string
\begin{picture}(100,10)(0,3)
\put(20,5){\circle*{4}} \put(40,5){\circle*{4}}
\put(80,5){\circle*{4}} \put(10,5){\line(1,0){40}}
\put(70,5){\line(1,0){20}} \put(62,5){\makebox(0,0){$\cdots$}}
\end{picture}
. The number of vertices of degree three is exactly $g$. In
general, any  vertex $v$ is decorated by the self-intersection of
the corresponding irreducible exceptional divisor $E_v$. In the
above diagram we put only the decoration of the vertex $v_0$,
which corresponds to the {\em unique} $(-1)$-curve, and which
also supports the strict transform of $\{f=0\}$. The other
decorations are not really essential in our next discussions; for
the description of the complete graph see e.g. \cite{BK} or
\cite{DSI}, section 4.I.

The above diagram can also be identified with the plumbing graph
of $(S^3,K_f)$. In this case the self-intersections are the
corresponding Euler numbers of the $S^1$-bundles, all the
surfaces used in the plumbing construction are $S^2$'s, and $K_f$
is a generic fiber of the unique $(-1)$-bundle.

The above resolution graph $G(f)$ will be denoted by the
following schematic diagram:

\begin{picture}(400,50)(0,0)
\put(250,25){\circle*{4}} \put(310,25){\makebox(0,0){$K_f$}}
\put(165,5){\framebox(100,40){}} \put(250,25){\vector(1,0){30}}
\put(250,35){\makebox(0,0){$-1$}}
\put(250,15){\makebox(0,0){$v_0$}}
\put(180,25){\makebox(0,0){$\bar{G}$}}
\dashline{3}(250,25)(235,30) \dashline{3}(250,25)(235,20)
\put(100,25){\makebox(0,0){$G(f):$}}
\end{picture}

The polynomial $\Delta(t)$ can also be deduced from $G(f)$ by
A'Campo's formula \cite{AC}. Notice that (if we disregard the
arrowhead) the graph $\bar{G}$ can be blown down completely.

We emphasize again, the information codified in the following
objects --  isotopy class of  $K_f\subset S^3$, set of Newton
pairs, set of linking pairs, Alexander polynomial, semigroup or
the embedded resolution graph -- are completely equivalent. This
means that the polynomial $Q$ introduced below can be deduced
from any of them.

\subsection{Definition of the polynomial $Q$}\label{a2} One has
the following identity  connecting the Alexander polynomial
$\Delta$ and the semigroup $\Gamma$,  cf. \cite{GZDC}:
$$\frac{\Delta(t)}{1-t}=\sum_{k\in \Gamma}\ t^k.$$
Since $\Delta(1)=1$ and $\Delta'(1)=\delta$ (use e.g. the formula
of  \ref{ap}), one gets
$$\Delta(t)=1+\delta(t-1)+(t-1)^2\cdot Q(t)$$
for some polynomial $Q(t)=\sum _{i=0}^{\mu-2} \alpha_it^i$ of
degree $\mu-2$ with integral coefficients. In fact, all the
coefficients $\{\alpha_i\}_{i=0}^{\mu-2}$ are strict positive,
and:
$$\delta=\alpha_0\geq \alpha_1\geq \cdots \geq \alpha_{\mu-2}=1.$$
Indeed, by the above identity, one has
$\delta+(t-1)Q(t)=\sum_{k\not\in\Gamma}t^k$, or
$Q(t)=\sum_{k\not\in\Gamma}(t^{k-1}+\cdots +t+1)$. This shows that
$$\alpha_i=\#\{k\not\in\Gamma\ :\ k>i\}.$$

\subsection{The surgery 3-manifold $M=S^3_{-p/q}(K)$}\label{a3}
In the sequel we fix an oriented algebraic knot $K_f\subset S^3$
and we consider  the oriented manifold $M:=S^3_{-p/q}(K_f)$
obtained by $-p/q$-surgery along $K_f\subset S^3$, where $p/q>0$
is a positive rational number ($p>0$, gcd$(p,q)=1$).

It is easy to see (e.g. by Kirby calculus, see e.g. \cite{GS})
that $M$ can be represented by the following plumbing diagram
(the symbols $v_0,\ldots , v_s$ are the `names' of the
corresponding vertices, they are not really parts of the
decoration of the diagram):

\begin{picture}(400,50)(120,0)
\put(250,25){\circle*{4}}  \put(165,5){\framebox(100,40){}}
\put(250,25){\line(1,0){50}} \put(300,25){\circle*{4}}
\put(300,15){\makebox(0,0){$v_1$}}
\put(350,15){\makebox(0,0){$v_2$}}

 \put(450,15){\makebox(0,0){$v_{s-1}$}}
 \put(500,15){\makebox(0,0){$v_s$}}

\put(250,35){\makebox(0,0){$-1$}}
\put(250,15){\makebox(0,0){$v_0$}}
\put(180,25){\makebox(0,0){$\bar{G}$}}
\put(140,25){\makebox(0,0){$G(M):$}} \dashline{3}(250,25)(235,30)
\dashline{3}(250,25)(235,20)

\put(300,25){\circle*{4}} \put(350,25){\circle*{4}}
\put(450,25){\circle*{4}} \put(500,25){\circle*{4}}
\put(300,25){\line(1,0){60}} \put(500,25){\line(-1,0){60}}
\put(405,25){\makebox(0,0){$\cdots$}}
\put(300,35){\makebox(0,0){$-k_1-m_f$}}
\put(350,35){\makebox(0,0){$-k_2$}}
\put(500,35){\makebox(0,0){$-k_s$}}
\put(450,35){\makebox(0,0){$-k_{s-1}$}}

\end{picture}

\noindent  where $k_1\geq 1$ and $k_j\geq 2$ $(2\leq j\leq s)$
are integers determined by the continued fraction
$$\frac{p}{q}=[k_1,k_2,\ldots,k_s]:=k_1-\cfrac{1}{k_2-\cfrac{1}{\ddots -\cfrac{1}{k_s}}}\ ;$$
and
$$m_f=a_gp_g,$$
or equivalently, in the language of the embedded resolution graph
of $f$ (cf. \ref{resgr}), $m_f$ is the multiplicity of the pull
back of the germ $f$ along the $(-1)$-curve. Topologically,
 e.g. if one uses Kirby calculus and blows down all the vertices
of $\bar{G}$, $m_f$ is a sum $\sum_{j}m_j^2$ of squares of a
sequence of linking numbers $m_j$. In fact, in the language of
plane curve singularities, this sequence of $m_j$'s is exactly
the `multiplicity sequence' of $f$.

Notice that if we would start with the unknot $K_f\subset S^3$,
then $M=S^3_{-p/q}(K_f)$ would be the lens space $L(p,q)$ (in
general, normalized with $0<q<p$, hence $k_1\geq 2$ as well).
This case is completely solved, see e.g. \cite{OSZINV}, section
10. Therefore, in the sequel we will assume that $K_f$ is the
unknot; and we accept the case $0<p/q\leq 1$  (i.e. $k_1=1$) as
well.

Take a point  $*\in S^3\setminus K_f$ and identify $S^3\setminus
*$ with $\R^3$. Then a Kirby diagram of $M$ is given by the knot
$K_f\subset \R^3$ with decoration (surgery coefficient) $-p/q$.
In particular, the Kirby diagram of the manifold $-M$ ($M$ with
reversed orientation) -- whose Heegaard Floer homology will be
computed in the sequel -- is given by the mirror image $K_f^m$ of
$K_f$ (with respect to any plane in $\R^3$)  with surgery
coefficient $p/q$; i.e. $-M=S^3_{p/q}(K_f^m)$.

\subsection{Some properties of the graph
$G(M)$.}\label{peculiar}

\subsubsection{}\label{h1} The plumbing graph $G(M)$ is a negative definite
connected tree.  For some particular plumbed 3-manifolds, with
such a graph, the Heegaard Floer homology is determined in
\cite{OSzP} and \cite{OSZINV}. Notice that, in general, $G(M)$
may have many `bad vertices' (in fact, their number is  $\#\{1\leq
i\leq s\,:\, p_i<q_i\}$), hence \cite{OSzP} cannot be applied.
Nevertheless, \cite{OSZINV} works. Indeed, first recall (cf.
\cite{OSZINV}) that a graph is called {\em almost rational},  or
$AR$, if it has at least one vertex with the following property:
if one replaces the decoration (Euler-number) of the
corresponding vertex by a smaller integer, then one gets a
rational graph (in the sense of Artin \cite{Artin62}). Then the
Heegaard Floer homology of a plumbed manifolds associated with
$AR$-graphs can be determined by \cite{OSZINV}.

\subsubsection{}\label{lemmaAR} {\bf Lemma.} {\em $G(M)$ is an
$AR$-graph.}

\begin{proof} There exists  a particular family of rational
graphs -- the so called {\em sandwiched graphs} (associated with
sandwiched singularities) -- which are characterized as follows:
they are subgraphs of (not necessarily minimal) connected negative
definite plumbing  graphs which can be completely blown down (cf.
e.g. \cite{Spiv}). We will show that if we modify the $-1$
decoration of $v_0$ in $G(M)$ into $-2$ (let us call this new
graph by $G(M)_{-2}$), we get a sandwiched graph. Indeed,
consider the graph $\bar{G}$. Blow up the $-1$ vertex $v_0$. The
new decoration of $v_0$ will be $-2$, while a new $-1$ vertex is
created. Then blow up this new vertex $(m_f+k_1-1)$ times. Then
its new decoration will be $-m_f-k_1$, while it has $(m_f+k_1-1)$
neighbours which are all $-1$ curves. Fix one of them, and blow up
$k_2-1$ times. If one continues this procedure, one  gets a graph
which contains $G(M)_{-2}$ as a subgraph.
\end{proof}

Using the results of \cite{OSZINV}, the above lemma has also the
following consequence: the plumbed 3-manifold associated with
$G(M)_{-2}$ is an $L$-space (in the sense of Ozsv\'ath-Szab\'o,
i.e. its reduces Heegaard Floer homology is trivial).


\section{Graded roots}

\subsection{Preliminary remark}\label{pr} The Heegaard Floer
homology of $M$ (in fact, of $-M$) will be computed in a
combinatorial way from the graph $G(M)$. Following \cite{OSZINV},
we do this via some intermediate objects, the $d$ graded roots
associated with $G(M)$, each corresponding to a
$spin^c$-structure of $M$. As we will see, the graded roots
contain more information than the Heegaard Floer homology, a fact
which hopefully will be more exploited in the future. On the
other hand, compared with the homology, they have an enormous
advantage: they can  much-much easier be  identified and described
(because, probably, they preserve the geometric information of
$M$ at more canonical level).

In this section we provide a short review of abstract graded
roots; for more details, see \cite{OSZINV}.

\subsection{Definition of a graded root $(R,\chi)$.}\label{2.1} \ \\
(1)  Let $R$ be an infinite tree with vertices $\calv$ and edges
$\cale$. We denote by $[u,v]$ the edge with
 end-points $u$ and $v$.  We say that $R$ is a graded root
with grading $\chi:\calv\to \Z$ if

(a) $\chi(u)-\chi(v)=\pm 1$ for any $[u,v]\in \cale$;

(b) $\chi(u)>\min\{\chi(v),\chi(w)\}$ for any $[u,v],\
[u,w]\in\cale$;

(c) $\chi$ is bounded below, $\chi^{-1}(k)$ is finite for any
$k\in\Z$, and $\#\chi^{-1}(k)=1$ if $k$ is sufficiently large.

\vspace{1mm}

\noindent (2) $v\in\calv $ is a {\em local minimum point} of the
graded root $(R,\chi)$ if $\chi(v)<\chi(w)$ for any edge $[v,w]$.
Their set is denoted by $\calv_{lm}$.  In fact, $\calv_{lm}$
coincides with the set of vertices with adjacent degree one.

\vspace{1mm}

\noindent  (3) A geodesic path connecting two vertices is {\em
monotone } if $\chi$ restricted  to the set of vertices on the
path is strict monotone. If a vertex $v$ can be connected by
another vertex $w$ by a monotone geodesic and $\chi(v)>\chi(w)$,
then we write $v\succ w$. $\succ$ is an ordering of $\calv$. For
any pair $v,w\in\calv$ there is a unique $\succ$-minimal vertex
$\sup(v,w)$ which dominates both.

\subsection{Examples.}\label{2.3}

(1) For any integer $n\in\Z$, let  $R_n$ be the tree with
$\calv=\{v^{k} \}_{ k\geq n}$ and
$\cale=\{[v^{k},v^{k+1}]\}_{k\geq n}$. The grading is
$\chi(v^{k})=k$.

(2) Let $I$ be a finite index set. For each $i\in I$ fix  an
integer $n_i\in \Z$; and for each pair $i,j\in I$ fix
$n_{ij}=n_{ji}\in\Z$ with the next properties: (i) $n_{ii}=n_i$;
(ii) $n_{ij}\geq \max\{n_i,n_j\}$; and (iii) $n_{jk}\leq
\max\{n_{ij},n_{ik}\}$ for any $ i,j,k\in I$.

For any $i\in I$ consider $R_i:=R_{n_i}$ with vertices
$\{v_i^{k}\}$ and edges $\{[v_i^{k},v_i^{k+1}]\}$, $(k\geq n_i)$.
In the disjoint union $\coprod_iR_i$,  for any pair $(i,j)$,
 identify $v_i^{k}$ and $v_j^{k}$,
resp. $[v_i^{k},v_i^{k+1}]$  and $[v_j^{k},v_j^{k+1}]$, whenever
$k\geq n_{ij}$. Write $\bar{v}_i^{k}$ for the class of $v_i^k$.
Then $R:=\coprod_iR_i/_\sim$ is a graded root with
$\chi(\bar{v}_i^{k})=k$.

Clearly $\calv_{lm}$ is a subset of $\{\bar{v}_i^{n_i}\}_{i\in
I}$. They are equal if in (ii) all the inequalities are strict.
(Otherwise some indices are superfluous, and the corresponding
$R_i$'s produce no contribution.)

(3) Any map $\tau:\{0, 1,\ldots,T_0\}\to \Z$ produces a starting
data for construction (2). Indeed, set $I=\{0,\ldots,T_0\}$,
$n_i:=\tau(i)$ ($i\in I$), and $n_{ij}:=\max\{n_k\,:\, i\leq
k\leq j\}$ for $i\leq j$. Then  $\coprod_iR_i/_\sim $ constructed
in (2) using this data will be denoted by $(R_\tau,\chi_\tau)$.

\vspace{2mm}

To any graded root, one can associate  a natural $\Z[U]$-module.

\subsection{Notations.}\label{zu}
Consider the $\Z[U]$-module $\Z[U,U^{-1}]$, and (following
\cite{OSzP}) denote by
 $\calt_0^+$ its quotient by the submodule  $U\cdot \Z[U]$.
It is a $\Z[U]$-module with grading  $\deg(U^{-h})=2h$.
Similarly, for any $n\geq 1$, define the $\Z[U]$-module
$\calt_0(n)$ as the quotient of $\Z\langle U^{-(n-1)},
U^{-(n-2)},\ldots, 1,U,\ldots \rangle$ by $U\cdot \Z[U]$ (with
the same grading). Hence, $\calt_0(n)$, as a $\Z$-module, is the
free $\Z$-module $\Z \langle 1,U^{-1},\ldots,U^{-(n-1)} \rangle$
(generated by $1,U^{-1},\ldots,U^{-(n-1)}$), and has finite
$\Z$-rank $n$.

More generally,  fix an arbitrary $\Q$-graded $\Z[U]$-module $P$
with $h$-homogeneous elements $P_h$. Then for any $r\in \Q$  we
denote by $P[r]$ the same module graded in such a way that
$P[r]_{h+r}=P_{h}$. Then set $\calt^+_r:=\calt^+_0[r]$ and
$\calt_r(n):=\calt_0(n)[r]$.

\subsection{Definition. The $\Z[U]$-module associated with a graded root.}\label{2.6}
For any graded root $(R,\chi)$, let $\bH(R,\chi)$ be the set of
functions $\phi:\calv\to \calt^+_0$ with the following property:
whenever $[v,w]\in \cale$ with $\chi(v)<\chi(w)$, then
\begin{equation*}
U\cdot \phi(v)=\phi(w).
\end{equation*}
Then $\bH(R)$ is a $\Z[U]$-module via $(U\phi)(v)=U\cdot \phi(v)$.
Moreover, $\bH(R)$ has a grading: $\phi\in \bH(R)$ is homogeneous
of degree $h\in\Z$ if for each $v\in\calv$ with $\phi(v)\not=0$,
$\phi(v)\in\calt^+_0$ is homogeneous  of degree $h-2\chi(v)$.

\vspace{2mm}

In fact, $\bH(R,\chi)$ can also be computed as follows (cf. (3.6)
in \cite{OSZINV}):

\subsection{Proposition.}\label{2.7} {\em Let $(R,\chi)$ be a
graded root. We order $\calv_{lm}$ as follows. The first element
$v_1$ is an arbitrary vertex with $\chi(v_1)=\min_{v}\chi(v)$. If
$v_1,\ldots, v_k$ is already determined, and
$J_k:=\{v_1,\ldots,v_k\}\varsubsetneq \calv_{lm}$, then let
$v_{k+1}$   be an arbitrary vertex in $\calv_{lm}\setminus J_k$
with $\chi(v_{k+1})=\min _{v\in \calv_{lm}\setminus J_k}\chi(v)$.
Let $w_{k+1}\in \calv$ be the unique $\succ$-minimal vertex of
$R$ which dominates both $v_{k+1}$, and at least one vertex from
$J_k$.  Then one has the following isomorphism of $\Z[U]$-modules
$$\bH(R,\chi)=\calt^+_{2\chi(v_1)}\oplus \oplus_{k\geq 2}\calt_{2\chi(v_k)}\big(
\chi(w_k)-\chi(v_k)\big).$$ In particular, with the notations
$m:=\min_v\chi(v)$ and $\bH_{red}(R,\chi):=\oplus_{k\geq
2}\calt_{2\chi(v_k)}\big( \chi(w_k)-\chi(v_k)\big)$, one has a
direct sum decomposition of graded $\Z[U]$-modules:}
$$\bH(R,\chi)=\calt^+_{2m}\oplus \bH_{red}(R,\chi).$$

\vspace{2mm}

Although the above proposition is elementary and in any concrete
example is very easy to apply, in some general situation can be
rather `inconvenient' to write down the summands. E.g., in the
case of \ref{2.3}(3), it is more natural to provide the integers
$\tau(i)$ (which eventually  have some geometric meaning) than to run
the above algorithm \ref{2.7}.  Therefore, sometimes we prefer to
stop at the level of $(R,\chi)$, but the interested reader can
always complete the picture with \ref{2.7}.

Proposition \ref{2.7} for  $(R_\tau,\chi_\tau)$ gives the
following:

\subsection{Corollary.}\label{2.10}  (See (3.8) in \cite{OSZINV}.) {\em Let
$(R_\tau,\chi_\tau)$  be a graded root associated with a function
$\tau:\{0,\ldots,T_0\}\to \Z$, cf. \ref{2.3}(3), which satisfies
$\tau(1)>\tau(0)=0$.  Then the $\Z$-rank of \ $\bH_{red}(
R_\tau,\chi_\tau)$ is:
$$\rank_\Z \bH_{red}(R_\tau)=\min_{i\geq 0}\tau(i)+\sum_{i\geq 0}\,
\max\{ \tau(i)-\tau(i+1),0\}.$$ Moreover, the summand
$\calt^+_{2m}$ of \ $\bH(T_\tau,\chi_\tau)$ has index
$m=\min_{i\geq 0}\tau(i)=\min_v\chi_\tau(v)$.}

\section{The Heegaard Floer homology $HF^+(-M)$}

\subsection{Preliminaries.}\label{invM} \
For any oriented rational homology 3-sphere $M$ the Heegaard Floer
homology $HF^+(M)$ was introduced by  Ozsv\'ath and Szab\'o in
\cite{OSz} (cf. also with their long list of articles). In fact,
$HF^+(M)$ is a $\Z[U]$-module with a $\Q$-grading compatible with
the $\Z[U]$-action, where $\deg(U)=-2$. Additionally, $HF^+(M)$
also has an (absolute) $\Z_2$-grading; $HF^+_{even}(M)$,
respectively $HF^+_{odd}(M)$, denote the part of $HF^+(M)$ with
the corresponding parity.

Moreover, $HF^+(M)$  has a natural direct sum decomposition of
$\Z[U]$-modules (compatible with all the gradings) corresponding
to the $spin^c$-structures of $M$:
$$HF^+(M)=\oplus_{\sigma\in Spin^c(M)}\ HF^+(M,\sigma).$$
For any $spin^c$-structure $\sigma$, one has a graded
$\Z[U]$-module isomorphism
$$HF^+(M,\sigma)=\calt^+_{d(M,\sigma)}\oplus HF^+_{red}(M,\sigma),$$
where $HF^+_{red}(M,\sigma)$ has a finite $\Z$-rank and an induced
(absolute) $\Z_2$-grading (and $d(M,\sigma)$ can also be defined
via this isomorphism). One also considers
$$\chi(HF^+(M,\sigma)):=\rank_\Z HF^+_{red,even}(M,\sigma)
-\rank_\Z HF^+_{red,odd}(M,\sigma).$$ Then one recovers the
(modified) Seiberg-Witten topological invariant of $(M,\sigma)$
(see \cite{Rus}) via
$$\ssw^{OSz}(M,\sigma):=\chi(HF^+(M,\sigma))-\frac{d(M,\sigma)}{2}.$$
With respect to the change of orientation the above invariants
behave as follows:
$$d(M,\sigma)=-d(-M,\sigma) \ \  \mbox{and} \ \
\chi(HF^+(M,\sigma))=-\chi(HF^+(-M,\sigma)).$$ Here one has a
natural identification of the $spin^c$ structures of $M$ and $-M$.
Notice also that one can recover $HF^+(M,\sigma)$ from
$HF^+(-M,\sigma)$ via (7.3) \cite{OSz} and (1.1) \cite{OSzTr}.

\subsection{}\label{det} Assume now that $M=S^3_{-p/q}(K_f)$ as
above. Then
 $H:=H_1(M,\Z)=\Z_p$. There is natural identification of
the $spin^c$-structures of $M$ with the classes of $\Z_p$  -- or,
with the  integers $a=0,1,\ldots, p-1$, $a=0$ corresponding to
the `canonical' $spin^c$-structure (cf. \cite{SWI}).
For the precise identification, see \ref{charel} and \ref{631}.
Let $\sigma_a$ be the $spin^c$-structure associated with $a$.

  Next theorem provides a
purely combinatorial description of $HF^+(-M)$ from $G(M)$. More
precisely, $HF^+(-M,\sigma_a)$ will be expressed in terms of the
delta-invariant $\delta$, the coefficients $\{\alpha_i\}_i$ of the
polynomial $Q$ (cf. \ref{a2}), and the integers $p,q$ and $a$.
The symbol  $\bms(q,p) $ denotes the Dedekind sum (see e.g.
\ref{adjf}). The integer $q'$ is  uniquely determined by $1\leq
q'\leq p$ and $qq'\equiv 1\ (mod\ p)$.

\subsection{Theorem.}\label{Maintheor} {\em For each fixed $a=0,1,\ldots,
p-1$, -- corresponding to the $p$ different $spin^c$-structures of
$M$ -- one defines the following objects\,:

\vspace{1mm}

$\bullet$ \ \ $t_a:=\Big[
\frac{(2\delta-1)q-a-1}{p}\Big];$

\vspace{1mm}

 $\bullet$ \ \  $r_a:= 3\bms(q,p) +2\sum
_{j=1}^a\Big\{\frac{jq'}{p}\Big\} -\frac{(1+2a)(p-1)}{2p}+\delta
\Big(1-\frac{q+1}{p}\Big)+\frac{\delta^2 q}{p}-\frac{2\delta
a}{p}$;

\vspace{1mm}

$\bullet$ \ \ a function $\tau_a:\{0,1,\ldots, 2t_a+2\}\to \Z$ by
}
$$\left\{ \begin{array}{l}
\tau_a(2t) =t(1-\delta)+\sum_{j=0}^{t-1}\Big[\frac{jp+a}{q}\Big], \ \ \ (t=0,\ldots, t_a+1);\\ \ \\
\tau_a(2t+1) = \tau_a(2t+2)+\alpha_{[(tp+a)/q]}, \ \  \
(t=0,\ldots, t_a).\end{array}\right.$$

\vspace{1mm}

{\em  $\bullet$ and the graded root $(R_{\tau_a},\chi_{\tau_a})$
associated with $\tau_a$.

\vspace{1mm}

 \noindent Then the following facts hold:

\vspace{1mm}

\ \ \ $HF^+_{odd}(-M,\sigma_a)=0$;

\vspace{1mm}

\ \ \ $HF^+_{even}(-M,\sigma_a)=\bH
(R_{\tau_a},\chi_{\tau_a})[r_a]$;

\vspace{1mm}

\ \ \ $d(-M,\sigma_a)=2\cdot \min \tau_a+r_a$.}

\vspace{2mm}

The proof is given in sections 5-6-7. We continue with some
remarks and corollaries.

\subsection{Remark.}\label{rema} Notice that for any
$t\in\{0,\ldots, t_a\}$, $\alpha_{[(tp+a)/q]}$ is strict positive,
hence $\tau_a(2t+1)>\tau_a(2t+2)$. On the other hand, using
properties of $\Gamma$ (see e.g. \ref{chp}, or \ref{taua}(2)), one
can verify that the following identity also holds:
$$\tau_a(2t+1)=\tau_a(2t)+\#\{\gamma\in \Gamma\,:\, \gamma\leq \,
(tp+a)/q\,\}.$$ Therefore, since $0\in \Gamma$, one has
$\tau_a(2t+1)>\tau_a(2t)$. In particular, the above representation
of the graded root is the most `economical': all the values are
essential, cf. also with \ref{2.3}(2).

This also shows that $(R_{\tau_a},\chi_{\tau_a})$ has exactly
$t_a+2$ local minimum points, and they correspond to the values
$\tau_a(2t), \ t=0,1,\ldots, t_a+1$.

\vspace{2mm}

\ref{Maintheor}, \ref{rema}  and \ref{2.10} imply the next two
corollaries.

\subsection{Corollary.}\label{sw} \
$$\ssw^{OSz}(M,\sigma_a)=-\ssw^{OSz}(-M,\sigma_a)=\frac{r_a}{2}
-\sum_{t\geq 0}\, \alpha_{[(tp+a)/q]}.$$

\vspace{2mm}

Recall that $HF^+(-M,\sigma_a)$ is a $\Z[U]$-module. Denote by
$\Ker U (\sigma_a)$ (resp. by $\coker\, U(\sigma_a)$) the kernel
(resp. the cokernel) of the $U$-action.  They are finitely
generated graded free $\Z$-module. Let $\Z_{(r)}$ denote a rank
one free $\Z$-module whose grading is $r$.

\subsection{Corollary.}\label{kerU} \
$$\Ker U(\sigma_a)=\bigoplus _{0\leq t\leq t_a+1}\ \Z_{(2\tau_a(2t)+r_a)}.$$
{\em In particular, $\Ker U(\sigma_a)$ depends only on the
integers $p$, $q$ $a$, and $\delta$, but not on the coefficients
$\{\alpha_i\}_i$. On the other hand,
$$\coker \, U(\sigma_a)=\bigoplus _{0\leq t \leq t_a}\
\Z_{(2\tau_a(2t+1)+r_a-2)},$$ which depends essentially on the
coefficients $\{\alpha_i\}_i$ of $Q$. Moreover:}
$$\rank_\Z\Ker U=\sum_{a=0}^{p-1} (t_a+2)=p+(2\delta-1)q=p+
\rank_\Z\coker U.$$

\subsection{Remark.}\label{exam}  If  $a\geq (2\delta-1)q$, then
$t_a=-1$, hence $R_{\tau_a}=R_0$ (cf. \ref{2.3}(1)). In
particular, for such $a$, $\bH_{red}(R_{\tau_a})=0$. For all the
other $a$'s,  $\bH_{red}(R_{\tau_a})$ is not trivial. E.g., for
the `canonical' $spin^c$-structure corresponding to $a=0$,
$\bH_{red}(R_{\tau_0})$ is never trivial. In particular, $-M$ is
never an $L$-space.

\section{Review of \cite{OSZINV} }

\subsection{Notations of the plumbing representation.}\label{3.1}
Let $G$ be a connected negative definite plumbing $AR$-graph
which determines the plumbed 3-manifold $M$ (cf. \ref{h1}).
Write  $\calj$ for the index set of vertices of $G$.
Each vertex $j\in\calj$  has its
decoration $e_j$, the corresponding Euler number.
 Let $L$ be the
free $\Z$-module of rank $\#\calj$ with fixed basis
$\{b_j\}_{j\in\calj}$. On $L$ there is a natural bilinear form
$B=(b_i,b_j)_{i,j\in\calj}$, defined by $(b_j,b_j)=e_j$; and for
$i\not=j$ one takes
 $(b_i,b_j)=1$ if the corresponding vertices $i,j$ are connected
by an edge in  $G(M)$, otherwise $=0$.
Sometimes we write $x^2$ for $(x,x)$.

For a cycle $x=\sum_j  n_jb_j\in L$ we write $x\geq 0$, if
$n_j\geq 0$ for any $j$.  Similarly,  $x\geq y$ if $x-y\geq 0$.

Set $L'=\Hom_\Z(L,\Z)$, and consider the exact sequence $0\to
L\stackrel{i}{\to} L'\to H\to 0$, where $i(x)=(x,\cdot )$. Notice
that $H=H_1(M,\Z)$. We will identify  the lattice $L'$ with a
sub-lattice of $L_\Q=L\otimes \Q$\,:   any $\alpha\in L'$ is
identified with $y_\alpha\in L_\Q$ which satisfies
$\alpha(x)=(y_\alpha , x)$ for any $x\in L$. Here
$(\cdot,\cdot)$  is  the natural extension of the form $B$  to
$L_\Q$.

Sometimes we will also use the dual base $\{g_j\}_{j\in \calj}$
of $L'$, i.e. $(g_i,b_j)=\delta_{ij}$ for any $i,j\in\calj$.
Clearly, $g_j$ is the $j^{th}$ column of $B^{-1}$ and
$(g_i,g_j)=B^{-1}_{ij}$. Since $B$ is negative definite, $-g_j\geq
0$ for any $j\in\calj$.

For any rational cycle $x=\sum_jr_jb_j$ we write $pr_j(x)$ for
the coefficient $r_j$.

\subsection{Characteristic elements.}\label{charel}
The set of characteristic elements are defined by
$$Char:=\{k\in L': \, k(x)+(x,x)\in 2\Z \ \mbox{for any $x\in L$}\}.$$
There is a natural action of $L$ on $Char$ by $x*k:=k+2i(x)$
whose orbits are of type $k+2i(L)$. Obviously, $H$ acts freely and
transitively on the set of orbits by
$[l']*(k+2i(L)):=k+2l'+2i(L)$. One has a natural identification
of the $spin^c$-structures $Spin^c(M)$ of $M$ with the set of
orbits of $Char$ modulo $2L$; and this identification is
compatible with the action of $H$ on both sets. Therefore, in the
sequel, we think about $Spin^c(M)$ by this identification, hence
any $spin^c$-structure
 of $M$ will be represented by an orbit  $[k]:=k+2i(L)\subset Char$.
(For more details, see e.g. \cite{SWI}.)

 The {\em canonical  characteristic element} $K\in
Char$ is defined by the `adjunction equations' $(K,b_j)=-e_j-2$
for any $j\in\calj$. Then $Char=K+2L'$, i.e. any $k\in Char$ has
the form $k=K+2l'$ for some $l\in L'$. Hence, the equivalence
class $[k]=K+2l'+2L $ can be identified with the sub-lattice
$l'+L$ of $L_\Q$. Next, we will choose a distinguished element of
$l'+L$. For this first consider:
$$S_\Q:=\{x\in L_\Q\, :\, (x,b_j)\leq 0 \ \mbox{for any }\ j\in
\calj\}.$$ Since $B$ is negative definite, $x\geq 0$ for any
$x\in S_\Q$. Then one has the following fact:

\vspace{1mm}

 {\em For any fixed $[k]=K+2(l'+L)$, the intersection
$(l'+L)\cap S_\Q$ in $L_\Q$ admits a unique minimal element
$l'_{[k]}$. (Here minimality is considered with respect to the
ordering $\leq $ in $L_\Q$.)}

\vspace{1mm}

For any fixed class ($spin^c$-structure) $[k]$ we fix the {\em
distinguished representative } $k_r:=K+2l'_{[k]}$ from $[k]$.
Since $l'_{[K]}=0$, the distinguished representative in $[K]$ is
$K$ itself.

\subsection{The graded roots $(R_k,\chi_k)$.}\label{mc1} For
any $k\in Char$ one defines $\chi_k:L\to\Z$ by
$$\chi_k(x):=-(k+x,x)/2.$$
 For any
$n\in \Z$, one constructs a finite 1-dimensional simplicial
complex $\bar{L}_{k,\leq n}$ as follows. Its 0-skeleton is
$L_{k,\leq n}:=\{x\in L:\, \chi_k(x)\leq n\}$. For each $x$ and
$j\in \calj$ with $x,x+b_j\in L_{k,\leq n}$, we consider a unique
1-simplex with endpoints at $x$ and $x+b_j$ (e.g., the segment
$[x,x+b_j]$ in $L\otimes \R$). We denote the set of connected
components of $\bar{L}_{k,\leq n}$ by $\pi_0(\bar{L}_{k,\leq
n})$. For any $v\in \pi_0(\bar{L}_{k,\leq n})$, let $C_v$ be the
corresponding connected component of $\bar{L}_{k,\leq n}$.

Next, we define the  graded root  $(R_k,\chi_k)$ as follows. The
vertices $\calv(R_k)$ are $\cup_{n\in \Z} \pi_0(\bar{L}_{k,\leq
n})$. The grading $\calv(R_k)\to\Z$, still denoted by $\chi_k$, is
$\chi_k|\pi_0(\bar{L}_{k,\leq n})=n$.

If $v_n\in \pi_0(\bar{L}_{k,\leq n})$, and $v_{n+1}\in
\pi_0(\bar{L}_{k,\leq n+1})$, and $C_{v_n}\subset C_{v_{n+1}}$,
then $[v_n,v_{n+1}]$ is an edge of $R_k$. All the edges
$\cale(R_k)$ of $R_k$ are obtained in this way.

Theorems (4.8) and (8.3) of \cite{OSZINV} say:

\subsection{Theorem.}\label{isohh} {\em Assume that $G$ is an $AR$-graph.
Then, for any $[k]\in Spin^c(M)$, $HF^+_{odd}(-M,[k])=0$ and
$$HF^+_{even}(-M,[k])=\bH (R_{k_r},\chi_{k_r})\Big[-\frac{k_r^2+\#\calj}{4}\Big].$$
In particular,}
$$d(-M,[k])=-\max_{k'\in[k]}\frac{(k')^2+\#\calj}{4}=-
\frac{k_r^2+\#\calj}{4}+2\min\chi_{k_r}.$$
We wish to emphasize that in most of the applications, the root
$(R_{k_r},\chi_{k_r})$ is {\em not} determined by its definition \ref{mc1},
but by  an algorithm which will be described in the sequel. (In
particular, the above definition, in concrete applications can be
neglected.)

\vspace{2mm}

Next, we  distinguish one of vertices of our $AR$-graph $G$,
which satisfies the definition of $AR$-graphs. By convention, it
corresponds to the index $0\in \calj$ (fixed for ever).

\subsection{The definition of the sequence $\{x_{[k]}(i)\}_{i\geq
0}$.}\label{xi} We fix a class $[k]$. Then for any integer $i\geq
0$ there exists a unique cycle $x_{[k]}(i)\in L$ with the
following properties:

(a) $pr_0(x_{[k]}(i))=i$;

(b)  $(x_{[k]}(i)+l'_{[k]},b_j)\leq 0$ for any $j\in
\calj\setminus \{0\}$;

(c) $x_{[k]}(i)$ is minimal (with respect the partial ordering
$\leq$) with the properties (a) and (b).

\subsection{Theorem.}\label{tauth} ((9.3) of \cite{OSZINV}) {\em
Fix a $spin^c$-structure $[k]$.  There exists an integer $T_0$
(which depend on $[k]$)  such that  $\chi_{k_r}(x_{[k]}(i+1))\geq
\chi_{k_r}(x_{[k]}(i))$ for any $i\geq T_0$. Moreover, the graded
root associated with $\tau_{[k]}:\{0,\ldots,T_0\}\to \N$ given by
$\tau_{[k]}(i):=\chi_{k_r}(x_{[k]}(i))$ satisfies}
$$(R_{k_r},\chi_{k_r})=(R_{\tau_{[k]}},\chi_{\tau_{[k]}}).$$

\subsection{Remark.}\label{inters} In general, it is not easy
to find the cycles $x_{[k]}(i)$.
Fortunately, one does not need all the coefficients of these
cycles, only the values $\chi_{k_r}(x_{[k]}(i))$. In most of the
cases they are computed inductively using the following equality
(cf. (9.1)(c) \cite{OSZINV}):
$$\chi_{k_r}(x_{[k]}(i+1))=\chi_{k_r}(x_{[k]}(i)+b_0).$$
Since the right hand side is
$\chi_{k_r}(x_{[k]}(i))+\chi_{k_r}(b_0)-(b_0,x_{[k]}(i))$,
basically one needs only the intersections $(b_0,x_{[k]}(i))$ for
any $i$.

\section{The first part  of the proof: $k_r^2+\#\calj$}

\subsection{The index set $\calj$} According to the shape of the
plumbing graph $G(M)$, the index set of its vertices is
$\calj=\bar{\calj}\cup\{1,\ldots,s\}$, where $\bar{\calj}$ is the
index set of the vertices of $\bar{G}$, and the indices
$\{1,\ldots, s\}$ correspond to $v_1,\ldots, v_s$. The
distinguish vertex of $G(M)$ is $v_0$ corresponding to
$0\in\bar{\calj}\subset \calj$. The base elements will also be
denoted accordingly: $b_0, b_1,\ldots, b_s$; or  $b_{j}$ for
$j\in\bar\calj$.

\subsection{The graph $_\delta\Gamma$.} As we already mentioned,
the sub-graph $\bar{G}$ can be
blown down. After this blow down, we get  a graph which will be
denoted as follows:

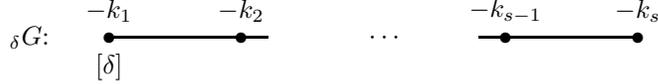
\begin{picture}(400,50)(250,0)
\put(300,25){\circle*{4}} \put(350,25){\circle*{4}}
\put(450,25){\circle*{4}} \put(500,25){\circle*{4}}
\put(300,25){\line(1,0){60}} \put(500,25){\line(-1,0){60}}
\put(405,25){\makebox(0,0){$\cdots$}}
\put(300,35){\makebox(0,0){$-k_1$}}
\put(350,35){\makebox(0,0){$-k_2$}}
\put(500,35){\makebox(0,0){$-k_s$}}
\put(450,35){\makebox(0,0){$-k_{s-1}$}}
\put(300,14){\makebox(0,0){$[\delta]$}}
\put(270,25){\makebox(0,0){$_{\delta}G$:}}
\end{picture}

We can think about this graph as the dual graph of a minimal
resolution (of a normal surface singularity) obtained from a
resolution with dual graph $G(M)$. In this sense, the new vertex
$v_1$ is a rational curve with a singular point which has
delta-invariant $\delta$; and the decorations $-k_i$ are the
corresponding self-intersections. On the other hand, one can
think about this graph also in the language of Kirby diagrams:
$v_1$ represents the knot $K_f\subset S^3$, the other vertices
represent unknots, they are linked as usual with linking number
one, and  $-k_i$ are the corresponding surgery coefficients.

The point is that a large number of numerical invariants of the
graph $G(M)$ -- needed for the descriptions of section 5 -- can
already be determined from $_{\delta}G$ in terms of $\delta$ and
$\{k_i\}_i$ (for this comparison, see the second part \ref{back} of this
section).
 The advantage of this is that the above graph $_{\delta} G$
is the same as the graph of a lens space $L(p,q)$, provided that
we disregard the decoration $[\delta]$. In particular, its
invariants can be computed by similar methods as those of lens
spaces -- in fact, in their computations we will even use the
corresponding relations valid for lens spaces. This will be the
subject of the first part of this section. Our model is section
10 of \cite{OSZINV} where we run the algorithm of section 5 for
lens spaces. The reader is invited to consult [loc. cit.] (and
also verify that our present claims, verified in [loc. cit.] for
case $q<p$, are valid for $q\geq p$ as well) .

In fact, any invariant  of the lattice (which does not involve
$\delta$) equals  the corresponding  invariant  of $L(p,q)$. But,
formulas which involve $\delta$ (like the `adjunction formula'
determining the canonical characteristic element, or the
`Riemann-Roch' formula for $\chi_k$) depend essentially on
$\delta$.

\subsubsection{}{\bf Notations.}  In order to make a distinction between
invariants of the graphs $G(M)$ and $_\delta G$, the invariants of
the second one will have an extra $\tilde{\ }$-decoration; e.g.
$\tilde{L}$ denotes the lattice of rank $s$ with base elements
$\{\tilde{b}_i\}_{i=1}^s$, $\{\tilde{g}_i\}_{i=1}^s\in\tilde{L}'$
are the dual bases, etc.

For any $1\leq i\leq j\leq s$ we write the continued fraction
$[k_i,\ldots, k_j]$ as a rational number $n_{ij}/d_{ij}$ with
$n_{ij}>0$ and $gcd(n_{ij},d_{ij})=1$. We also set $n_{i,i-1}:=1$
and $n_{ij}:=0$ for $j<i-1$. In fact, one gets $d_{ij}=n_{i+1,j}$.
 Clearly $n_{1s}=p$ and $n_{2s}=q$. We also
write $q':=n_{1,s-1}$. One can verify (e.g. by induction) that
$p\geq q'\geq 1$ (even if $p<q$), and $qq'\equiv 1 \ (mod\ p)$
(since $qq' =pn_{2,s-1}+1$).

\subsubsection{}\label{l3}{\bf The group $\tilde{H}$ and the elements
$\tilde{l}'_{[k]}$.}  Clearly
$\tilde{H}=\tilde{L}'/\tilde{L}=\Z_p$, and
$[\tilde{g}_s]=\tilde{g}_s+\tilde{L}$ is one of its generators.
Then $[\tilde{g}_j]=[n_{j+1,s}\tilde{g}_s]$ in $\tilde{H}$ $(1\leq
j\leq s)$. The set of $spin^c$-structures  is the set of orbits
$\{-a\tilde{g}_s+\tilde{L}\}_{0\leq a<p}$ (we prefer to use this
sign, since $-\tilde{g}_s>0$). For any $0\leq a<p$ write
$$\tilde{l}'_{[-a\tilde{g}_s]}=-(a_1\tilde{g}_1+a_2\tilde{g}_2+\cdots+a_s\tilde{g}_s).$$
 By \cite{OSZINV} \S 10, the system  $(a_1,\ldots, a_s)$ can be realized as the set of
coefficients of a minimal vector $\tilde{l}'_{[-a\tilde{g}_s]}$
for some $0\leq a<p$ if and only if   the entries satisfy the
inequalities:
$$(SI)\hspace{1cm} \left\{ \begin{array}{l}
a_1\geq 0, \ldots , a_s\geq 0\\
n_{i+1,s}a_i+n_{i+2,s}a_{i+1}+\cdots + n_{ss}a_{s-1}+a_s<n_{is}
 \ \mbox{for any $1\leq i\leq s$}.\end{array}\right.$$
The integer $0\leq a<p$ can be recovered from $(a_1,\ldots,a_s)$
(which satisfies (SI)) by
\begin{equation*}a=
n_{2s}a_1+n_{3s}a_2+\cdots + n_{ss}a_{s-1}+ a_s.\end{equation*}
Conversely, any $0\leq a<p$ determines inductively the entries
$a_1,\cdots , a_s$ by
$$a_i=\Big[ \frac{a-\sum_{t=1}^{i-1} n_{t+1,s}a_t}{n_{i+1,s}}\Big] \ \ (1\leq i\leq s).$$

\subsubsection{}\label{adjf}{\bf The canonical characteristic element.} Let
$_\delta\tilde{K}$ denote the  canonical characteristic element
associated with the graph $_\delta G$ (see below). It is
convenient to consider the canonical characteristic element
$\tilde{K}$ of the graph of $L(p,q)$ (which is obtained from
$_\delta G$ by omitting $\delta$). The adjunction relations
defining $\tilde{K}$ are the usual ones, but those which identify
$_\delta\tilde{K}$ are the following (see also \ref{pik}):
$(_\delta\tilde{K},\tilde{b}_j)-k_j+2$ should equal twice the
delta-invariant of the vertex $v_j$, namely $=2\delta$ for $j=1$,
and $=0$ otherwise. Hence:
$$_\delta\tilde{K}=\tilde{K}+2\delta\tilde{g}_1.$$
In particular,
$$_\delta\tilde{K}^2+s=\tilde{K}^2+s+4\delta(\tilde{K},\tilde{g}_1)
+4\delta^2(\tilde{g}_1,\tilde{g}_1).$$ $(\tilde{K},\tilde{g}_1)$
is the $\tilde{b}_1$-coefficient of $\tilde{K}$, which equals
(e.g. by \cite{SWI}, (5.2)) $-1 +(q+1)/p$. Similarly,
$\tilde{g}_1^2=\tilde{B}^{-1}_{11}=-n_{2s}/p=-q/p$.  Notice also
that (cf. (7.1) \cite{SWI}, or (10.4) of \cite{OSZINV} -- although
this formula goes back to the work of Hirzebruch):
$$\tilde{K}^2+s=\frac{2(p-1)}{p}-12\cdot \bms(q,p),$$
where $\bms(q,p)$ denotes the {\bf Dedekind sum}
\[
\bms(q,p)=\sum_{l=0}^{p-1}\Big(\Big( \frac{l}{p}\Big)\Big)
\Big(\Big( \frac{ ql }{p} \Big)\Big), \ \mbox{where} \ \
((x))=\left\{
\begin{array}{ccl}
\{x\} -1/2 & {\rm if} & x\in {\R}\setminus {\Z}\\
0 & {\rm if} & x\in {\Z}.
\end{array}
\right.
\]
Therefore,
$$_\delta\tilde{K}^2+s=\frac{2(p-1)}{p}-12\bms(q,p)-4\delta(1-\frac{q+1}{p})-4\delta^2\frac{q}{p}.$$
The distinguished characteristic element $\tilde{k}_r$ of the
orbit $-a\tilde{g}_s+\tilde{L}$ is
$$\tilde{k}_r=\ _\delta\tilde{K}+2\tilde{l}'_{[-a\tilde{g}_s]}.$$
Therefore,
$$\tilde{k}_r^2+s=\ _\delta\tilde{K}^2+s+4\cdot (\,
\tilde{K}+\tilde{l}'_{[-a\tilde{g}_s]},\tilde{l}'_{[-a\tilde{g}_s]}\,
)+8\delta \cdot (\,\tilde{g}_1, \tilde{l}'_{[-a\tilde{g}_s]}\,).$$
By \cite{OSZINV} (10.4.1) one has
$$(\,\tilde{K}+\tilde{l}'_{[-a\tilde{g}_s]},\tilde{l}'_{[-a\tilde{g}_s]}\,)=
\frac{a(p-1)}{p}-2\cdot \sum _{j=1}^a\Big\{\frac{jq'}{p}\Big\}.$$
On the other hand, by the above formulas: $$ (\,\tilde{g}_1,
\tilde{l}'_{[-a\tilde{g}_s]}\,)=-\sum_{i=1}^sa_i(\tilde{g}_1,\tilde{g}_i)=
\sum_{i=1}^s a_i\frac{n_{i+1,s}}{p}=\frac{a}{p}.$$ Summing all,
one gets
$$-\frac{\tilde{k}_r^2+s}{4}=
3\bms(q,p) +2\sum _{j=1}^a\Big\{\frac{jq'}{p}\Big\}
-\frac{(1+2a)(p-1)}{2p}+\delta
\Big(1-\frac{q+1}{p}\Big)+\frac{\delta^2 q}{p}-\frac{2\delta
a}{p}.$$

\subsection{Back to the graph $G(M)$.}\label{back} Finally, we compare the
invariants of the graphs $G(M)$ and $_\delta G$. Some of the
arguments are well-known, nevertheless for the convenience of the
reader we provide the proofs.

There are two natural morphisms connecting the corresponding
lattices $L$ and $\tilde{L}$. The first one, $\pi_*:L\to
\tilde{L}$ is defined by $\pi_*(b_j)=0$ for any
$j\in\bar{\calj}$, while $\pi_*(b_i)=\tilde{b}_i$ for
$i=1,\ldots, s$.

In order to define the second morphism, we need an additional
construction.

Let $Z_f:=\sum_{j\in\bar{\calj}}m_jb_j$ be the cycle supported by
$\bar{G}$ which satisfies
$$(Z_{f}+b_1,b_j)=0 \ \mbox{for $j\in\bar{\calj}$}.$$
Since $\det\bar{G}=\pm 1$, this system has a unique integral
solution $\{m_j\}_j$. In fact, in terms of the diagram $G(f)$ (cf.
\ref{alg1}), $m_j$ is the  vanishing order of the pull back of
$f$ along the corresponding irreducible exceptional divisor.
E.g., $m_0=m_f$.

Then one defines $\pi^*:\tilde{L}\to L$ by
$\pi^*(\tilde{b}_j)=b_j$ for $j\geq 2$ and $\pi^*(\tilde{b}_1)=
Z_f+b_1$. By the very definition follows the `projection formula':
$$(\pi^*(\tilde{l}),l)=(\tilde{l},\pi_*(l))\ \ \ \mbox{for} \ \tilde{l}\in\tilde{L},\ l\in L.$$

\subsubsection{}\label{631}{\bf The group $H$.} $\pi_*$ has a natural
extension to $L_\Q\to \tilde{L}_\Q$, and $\pi_*(L')\subset
\tilde{L}'$. Therefore, $\pi_*$ induces a group morphism $H\to
\tilde{H}$. Its surjectivity follows from $\pi_*\pi^*=1$. In
order to prove injectivity, consider an $l'\in L'$ with
$\pi_*(l')\in \tilde{L}$. Then $\pi^*\pi_*(l')$ is an element of
$L$. Since $\pi^*\pi_*(l')-l'$ is supported by $\bar{G}$, and
$\det\bar{G}=\pm 1$, one gets that $\pi^*\pi_*(l')-l'\in L$ as
well. Hence $l'\in L$, and $\pi_*$ induces an isomorphism $H\to
\tilde{H}=\Z_p$, and  $H$ is generated by $g_s$.

In the sequel, the set of $spin^c$-structures of $M$ will be
identified with the set of orbits $\{-ag_s+L\}_{0\leq a<p}$. We
denote by $\sigma_a$ that $spin^c$-structure which correspond to
the orbit $-ag_s+L$.

\subsubsection{}{\bf The element $l'_{[-ag_s]}$.} We claim that
$l'_{[-ag_s]}=\pi^*( \tilde{l}'_{[-a\tilde{g}_s]})$. Indeed, $
l':=\pi^*( \tilde{l}'_{[-a\tilde{g}_s]})\in S_\Q$ by the
projection formula. Next, we verify the minimality of $l'$.
Notice that any $l\in L$ with $l\geq 0$ can be written in the
form $\pi^*(\tilde{x})+\bar{x}$, where $\tilde{x}\geq 0$ and
$\bar{x}$ is supported by $\bar{G}$. Assume that for such
$\pi^*(\tilde{x})+\bar{x}\geq 0$ one has:
\begin{equation*}
(l'-\pi^*(\tilde{x})-\bar{x},b_j)\leq 0\ \mbox{ for any
$j\in\calj$}.\tag{$*$}\end{equation*} Since
$(\pi^*(\tilde{y}),b_j)=0$ for any $j\in\bar{\calj}$, one gets
that $(-\bar{x},b_j)\leq 0$ for any $j\in\bar{\calj}$. Since
$\bar{G}$ is negative definite, it follows that $\bar{x}\leq 0$.
Using this, ($*$) for $j=1,\ldots, s$ gives
$$0\geq (l'-\pi^*(\tilde{x}),b_j)+(-\bar{x},b_j)\geq
(l'-\pi^*(\tilde{x}),b_j)=(\tilde{l}'_{[-a\tilde{g}_s]}-\tilde{x},\tilde{b}_j).$$
hence, by the minimality of $\tilde{l}'_{[-a\tilde{g}_s]}$, one
has $\tilde{x}=0$. But then $\bar{x}\geq 0$ and $\bar{x}\leq 0$
implies $\bar{x}=0$ as well.

\subsubsection{}\label{pik}{\bf Claim: $\pi_*(K)=\  _\delta\tilde{K}$.} First
notice that by the projection formula:
$$\pi_*(g_j)=\left \{\begin{array}{ll}
m_j\tilde{g}_1 & \mbox{if $j\in\bar{\calj}$}\\
\tilde{g}_j & \mbox{if $j=1,\ldots, s$.}\end{array}\right.$$
(where $Z_f=\sum_{j\in \bar{\calj}}m_jb_j$ as above).
The adjunction relations for   $G(M)$ can be rewritten as
$$K=-\sum_{j\in \calj}b_j-\sum_{j\in \calj}(2-s_j)g_j,$$
where $s_j$ is the adjacent degree of the vertex $j$ in $G(M)$.
They also satisfy (cf. \cite{AC}):
$$\sum_{j\in\bar{\calj}}(2-s_j)m_j=1-\mu=1-2\delta.$$
Then, taking $\pi_*$ one gets
$\pi_*(K)=\tilde{K}+2\delta\tilde{g}_1=\ _\delta\tilde{K}$.


\subsubsection{}{\bf Claim: $K^2+\#\calj=\ _\delta\tilde{K}^2+s$.}
Indeed, set  $K_{\bar{G}}:=K-\pi^*(_\delta\tilde{K})$. Then
$K_{\bar{G}}$ is supported by $\bar{G}$ and satisfies the
adjunction relations $(K_{\bar{G}},b_j)=-e_j-2$ for any
$j\in\bar{\calj}$, hence it is the canonical cycle of $\bar{G}$.
Moreover, since $\bar{G}$ can be blown down (and since by an
elementary blow up of a smooth point $K^2 $ decreases by 1), one
has $K_{\bar{G}}^2=-(\#\calj-s)$. By projection formula $\pi^*(
_\delta\tilde{K})$ is orthogonal to $K_{\bar{G}}$, hence
$$K^2=(\pi^*( _\delta\tilde{K})+K_{\bar{G}})^2=\
_\delta\tilde{K}^2+K_{\bar{G}}^2=\
_\delta\tilde{K}^2-(\#\calj-s).$$ Finally, this claim, via the
projection formula, shows for any $spin^c$-structure that
$$k_r^2+\#\calj=\tilde{k}_r^2+s.$$

\subsection{Example. Integer surgery.}\label{int1} Assume
that $q=1$, i.e.  $M=S^3_{-p}(K_f)$. Then $q'=1$ as well, and
$$-\frac{k_r^2+\#\calj}{4}=-\frac{\tilde{k}_r^2+s}{4}=
\frac{(p+2\delta-2-2a)^2-p}{4p}.$$

\subsection{Example. $(1/q)$-surgery.}\label{1perq} Assume
that $p=1$, i.e.  $M=S^3_{-1/q}(K_f)$. Then again $q'=1$. There
is only one $spin^c$-structure corresponding to $a=0$, and
$$-\frac{K^2+\#\calj}{4}=-\frac{\tilde{K}^2+s}{4}=
q\delta(\delta -1).$$

\section{The second part of the proof:
$(R_{\tau_{[k]}},\chi_{\tau_{[k]}})$}

The construction of the cycles $x_{[k]}(i)$ is given in two steps.
The first step provides similar cycles associated with the graph
$\bar{G}$, and it is based on the combinatorial properties of the
 graph $G(f)$  (involving also some techniques of plane
curve singularities). In particular, we prefer to think about
$\bar{G}$ as the dual graph of irreducible exceptional curves
obtained by repeated blow ups of $(\C^2,0)$. We will use the
notation $S:=\{y:(y,b_j)\leq 0, \ j\in\bar{\calj}\}$.

\subsection{The cycles $\{y(i)\}_{i\geq 0}$.}\label{yi} Since  $\bar{G}$
 is negative definite, (7.6)
of \cite{OSZINV} guarantees, for any $i\geq 0$, the existence of a
positive cycle $y(i)\geq 0$ (supported by $\bar{G}$) with the
following properties:

(a) $pr_0(y(i))=i$;

(b)  $(y(i),b_j)\leq 0$ for any $j\in \bar{\calj}\setminus \{0\}$;

(c) $y(i)$ is minimal (with respect the partial ordering $\leq$)
with the properties (a) and (b).

\vspace{2mm}

 E.g., $y(0)=0$.
Although there is very precise algorithm which determines all the
cycles $y(i)$ (see e.g. the proof of \ref{lm1}, or \cite{OSZINV}),
we are not interested in all the coefficients of
$y(i)$. But what we really have to know is the set of intersection
numbers $(y(i),b_0)$ (cf. \ref{inters}).
  Let $Z_f$ be the divisorial part
(supported by $\bar{G}$) of the germ $f$, cf. \ref{back}, which
satisfies
$$(Z_f,b_j)=\left\{\begin{array}{rl}
0&\ \mbox{if} \  j\in \bar{\calj}\setminus \{0\};\\
-1& \ \mbox{if} \ j=0.\end{array}\right.$$ Recall also that
$\Gamma\subset\N$ denotes the semigroup of $f$, see \ref{alg1}.

\subsection{Proposition.}\label{propyi} {\em

(a) If $i=tm_f+i_0$ with $t\geq 0$ and $0\leq i_0<m_f$, then
$y(i)=tZ_f+y(i_0)$;

(b) For any $i<m_f$ one has}
$$(y(i),b_0)=\left\{\begin{array}{ll}
1&\mbox{if $i\not\in\Gamma$};\\
0&\mbox{if $i\in\Gamma$}.\end{array}\right.$$

\begin{proof} {\em (a)}
 Clearly, $y'(i):=tZ_f+y(i_0)$ satisfies (a)-(b), we need to
verify (c). But if for some $y\geq 0$,  $y'(i)-y$ satisfies
(a)-(b), then $y(i_0)-y$ would also satisfy (a)-(b) for $i_0$,
hence $y=0$ by minimality of $y(i_0)$. Part $(b)$ will follow from
the next sequence of lemmas.
\end{proof}

The first lemma was proved and used intensively in \cite{OSZINV}
as a general principle of rational graphs. For the convenience of
the reader, we sketch its proof, for more details, see [loc.
cit.].

\subsubsection{}{\bf Lemma.}\label{lm1} {\em $(y(i),b_0)\leq 1$ for any
$i\geq 0$.}

\begin{proof} We denote by $Z_{min}$ Artin's minimal
cycle associated with $\bar{G}$, i.e. the minimal cycle in $S$.
The cycle $Z_{min}$ can be determined by the following (Laufer's)
algorithm (\cite{Laufer72}): Construct a `computation' sequence of cycles
$\{Z_l\}_l$ as follows. Set $Z_0:=0$, $Z_1:=b_0$; if $Z_l$ is
already constructed, but for some $j\in\bar{\calj}$ one has
$(Z_l,b_j)>0$, then take $Z_{l+1}:=Z_l+b_j$. If $Z_l\in S$ then
stop and $Z_l=Z_{min}$.

Now, by a geometric genus computational algorithm, one gets for
all rational graphs that whenever $(Z_l,b_j)>0$ in the above
algorithm, in fact one has the equality ($\dagger$): $(Z_l,b_j)=1$
(cf. \cite{Laufer72}).
Since $\bar{G}$ can be blown down, it is rational, hence ($\dagger$) works.

The point is that $y(i)$ can also be computed by a similar
algorithm: Assume that $y(i)$ is already  determined. Then
construct the sequence of cycles $\{Z_l\}_l$ as follows.
$Z_0:=y(i)$, $Z_1:=y(i)+b_0$, and if (for $l\geq 1$) $(Z_l,b_j)>0$
for some $j\in\bar{\calj}\setminus \{0\}$, then take $Z_{l+1}:=
Z_l+b_j$, otherwise stop and write $Z_l=y(i+1)$.

Then one can verify that any sequence connecting $y(i)$ with
$y(i+1)$ can be considered as part of a computation sequence
associated with $Z_{min}$, hence  lemma follows by the above
property ($\dagger$).
\end{proof}

\subsubsection{}{\bf Lemma.}\label{lm2} {\em Fix an arbitrary $i\geq 0$.
Then $(y(i),b_0)\leq 0 $ if and only if \ $i\in\Gamma$.}

\begin{proof} `$\Rightarrow$'\ If $(y(i),b_0)\leq 0$ then
$y(i)\in S$. Since $\bar{G}$ is the dual graph of a modification
of $(\C^2,0)$, the cycle $y(i)$ is the divisorial part $Z_g$ of a
holomorphic germ $g\in {\mathcal O}_{(\C^2,0)}$. Since the
intersection multiplicity $i_0(f,g)$ of the germs $f$ and $g$ at
$0\in\C^2$ is the multiplicity of $Z_g$ along the exceptional
divisor supporting the arrow of $f$, one gets
$i_0(f,g)=pr_0(y(i))=i$. But, by definition, $i_0(f,g)\in\Gamma$.

`$\Leftarrow$'\ If $i\in \Gamma$, then there exists
$g\in{\mathcal O}_{(\C^2,0)}$ with $i_0(f,g)=i$; i.e. with
$pr_0(Z_g)=i$ (see above) and $Z_g\in S$ (a general property of
divisorial parts of holomorphic germs). Then the minimality of
$y(i)$ implies $y(i)\leq Z_g$. This, and $pr_0(y(i))=pr_0(Z_g)$,
and the fact that $B$ off the diagonal has positive entries imply
that $(b_0,y(i))\leq (b_0,Z_g)$. But, since $Z_g\in S$,
$(b_0,Z_g)\leq 0$.
\end{proof}

Lemma \ref{lm2} can be improved for $i<m_f$:

\subsubsection{}{\bf Lemma.}\label{lm3} {\em Assume that $i<m_f$. Then
$(y(i),b_0)=0$ if and only if \ $i\in \Gamma$.}

\begin{proof} Assume that $(y(i),b_0)<0$ for some $i$.
Then, by a verification $y(i)-Z_f\in S$, hence (since $\bar{G}$ is
negative definite) $y(i)-Z_f\geq 0$. In particular,
$i-m_f=pr_0(y(i)-Z_f)\geq 0$.
\end{proof}

\subsection{Back to $G(M)$.}\label{back2} We consider again the
graph $G(M)$. We fix an orbit $[k]=-ag_s+L$ for some $0\leq a<p$
(corresponding to the $spin^c$-structure $\sigma_a$ of $M$). In
the sequel we will use freely the notations of the previous
section. Additionally, we write for any $j=1,\ldots, s$:
$$a_j' :=n_{j+1,s}a_j+a_{j+2,s}a_{j+1}+\ldots+ a_s.$$
E.g., by \ref{l3}, $a'_1=a$.

 Consider now the
sub-graph of $G(M)$ with vertices $v_0,\ldots, v_s$ and the $s$
edges connecting them. We wish to identify the cycle $z(i)=ib_0
+\sum_{i=1}^su_jb_j$ which has the properties of  $x(i)$
`restricted on this sub-graph'. More precisely:

\subsection{Proposition.}\label{zi} {\em Fix $i\geq 0$ and $0\leq a<p$.

(a) For $j=1,\ldots, s$ define the integers $u_j$ inductively by:
$$u_1:=\Big\lceil\,\,\frac{iq-a}{p+qm_f}\,\,\Big\rceil,\ \
u_j:=\Big\lceil\,\,\frac{u_{j-1}n_{j+1,s}-a_j'}{n_{j,s}}\,\,\Big\rceil
\ (1<j\leq s).$$ Then the cycle $z(i):=ib_0+\sum_{j=1}^su_jb_j$ is
$\geq 0$ and  satisfies
$$(z(i)+l'_{[-ag_s]},b_j)\leq 0 \ \mbox{for any $j=1,\ldots ,s
$}.$$

(b) Assume that the cycle
$\bar{z}(i):=ib_0+\sum_{j=1}^s\bar{u}_jb_j$ is $\geq 0$ and
satisfies
$$(\bar{z}(i)+l'_{[-ag_s]},b_j)\leq 0 \ \mbox{for any $j=1,\ldots ,s
$}.$$ Then
$$\bar{u}_1\geq \Big\lceil\,\,\frac{iq-a}{p+qm_f}\,\,\Big\rceil\ ,\ \ \mbox{and }
\ \bar{u}_j\geq
\Big\lceil\,\,\frac{\bar{u}_{j-1}n_{j+1,s}-a_j'}{n_{j,s}}\,\,\Big\rceil
\ \mbox{for } \ 1<j\leq s.$$ }

\begin{proof} This property was used in a similar situation for
the plumbing graph of a Seifert manifolds, applied for its `legs',
cf. 11.11 of \cite{OSZINV}. For the convenience of the reader we
sketch the proof.

(a)  By the definition of $u_j$ one has $u_jn_{js}\geq
u_{j-1}n_{j+1,s}-a'_j$. This is equivalent, via the identities
$n_{js}=k_jn_{j+1,s}-n_{j+2,s}$ resp.
$a_j'=n_{j+1,s}a_j+a'_{j+1}$, with
$$u_jk_j-u_{j-1}+a_j\geq
\frac{u_jn_{j+2,s}-a'_{j+1}}{n_{j+1,s}}.$$ Using the definition
of $u_{j+1}$, one gets $u_jk_j-u_{j-1}+a_j\geq u_{j+1}$, which is
exactly $(z(i)+l'_{[-ag_s]},b_j)\leq 0 $. (b) can be proved by a
descending induction, using the above sequence of arguments in
reversed order.
\end{proof}

\subsection{Corollary.}\label{corxi} {\em Fix $0\leq a<p$ and
write $i=tm_f+i_0$ for some $t\geq 0$ and $0\leq i_0<m_f$. Then
$$x_{[k]}(i)=t\cdot Z_f +y(i_0)+\Big\lceil
\,\,\frac{iq-a}{qm_f+p}\,\,\Big\rceil\, b_1+\sum_{j\geq
2}u_jb_j.$$ In particular,
$$(x_{[k]}(i),b_0)=-t+\Big\lceil
\,\,\frac{iq-a}{qm_f+p}\,\,\Big\rceil\,+(y(i_0),b_0).$$ Moreover:}
$$\chi_{k_r}(\, x_{[k]}(i+1)\, )-\chi_{k_r}(\, x_{[k]}(i)\, )
=t+1-\Big\lceil \,\,\frac{iq-a}{qm_f+p}\,\,\Big\rceil\, -
\left\{\begin{array}{ll} 1&\mbox{if $i_0\not\in\Gamma$}\\
0&\mbox{if $i_0\in\Gamma$}.\end{array}\right.$$

\begin{proof} Since the subgraphs $\bar{G}$ and its complement are
connected in  $G(M)$ only through $v_0$, the two parts of
$x_{[k]}(i)$ -- one supported by $\bar{G}$, the other by its
complement -- have independent lives. But, for any
$j\in\bar{\calj}$, $(x_{[k]}(i)+l'_{[k]},b_j)=
(x_{[k]}+\pi^*(\tilde{l}'_{[k]}), b_j)=(x_{[k]},b_j)$. Hence
$x_{[k]}(i)$ restricted on $\bar{G}$ should be exactly $y(i)$.
Clearly, the restriction on the complement of $\bar{G}$ should be
$z_{[k]}(i)$. Hence the first two statements follow.  The last is
the consequence of \ref{inters} together with the identity
$\chi_{k_r}(b_0)=1$, for
$(l'_{[k]},b_0)=(\pi^*(\tilde{l}'_{[k]}),b_0)=0$.
\end{proof}

\subsection{$\tau_a$}\label{taua} Corresponding to this and \ref{tauth}, we write
$\tau_a(i):=\chi_{k_r}(x_{[k]}(i))$. With the notations of
\ref{corxi}, notice that
\begin{equation*}
\frac{iq-a}{qm_f+p}\leq t+1,\end{equation*} hence
$\tau_a(i+1)-\tau_a(i)\geq -1$ for any $i$, and $=-1$ only in
very special situations, namely if
\begin{equation*}
\frac{(tm_f+i_0)q-a}{qm_f+p}>t \ \ \mbox{and} \ \ \
i_0\not\in\Gamma.\tag{1}\end{equation*} In order to analyze when
is this possible, we will consider sequences $Seq(t):=\{tm_f+i_0:
0\leq i_0<m_f\}$ for fixed $t\geq 0$. In such a sequence, notice
that the very last element of $\N\setminus \Gamma$, namely $\mu-1=
2\delta-1$ is strict smaller than $m_f-1$ (a fact, which can be
proved by induction over the number of Newton pairs),
hence the complete set $\N\setminus \Gamma$
sits in $\{0, \ldots, m_f-1\}$. Therefore, there  exists in
$Seq(t)$ an $i_0$ satisfying (1) if and only if
\begin{equation*}
\frac{(tm_f+2\delta-1)q-a}{qm_f+p}>t.\end{equation*} This is
equivalent to $t\leq t_a$, for $t_a$ defined in \ref{Maintheor}.
In other words, if $i\geq T_0:=(t_0+1)m_f$, then $\tau_a(i+1)\geq
\tau_a(i)$, hence those values of $\tau_a$ provide no
contribution in the graded root (cf. \ref{2.3}). Moreover, for
$t\in\{0,\ldots,t_a\}$, in $Seq(t)$ one has:
\begin{equation*}
\Delta(i_0):=
\tau_a(tm_f+i_0+1)-\tau_a(tm_f+i_0)=\left\{\begin{array}{rll} 0&
\ \mbox{if} \ i_0\leq (tp+a)/q, & \ \mbox{and} \
i_0\not\in\Gamma;\\
 +1&
\ \mbox{if} \ i_0\leq (tp+a)/q, & \ \mbox{and} \
i_0\in\Gamma;\\
-1& \ \mbox{if} \ i_0 > (tp+a)/q, & \ \mbox{and} \
i_0\not\in\Gamma;\\
0& \ \mbox{if} \ i_0 > (tp+a)/q, & \ \mbox{and} \ i_0\in\Gamma.
\end{array}\right.\end{equation*}
In particular,  $\Delta(i_0)\geq 0$  for $0\leq i_0\leq (tp+a)/q$
with exactly $$A_t:=\#\{\gamma\in\Gamma: \gamma\leq (tp+a)/q\}$$
times taking the value $+1$, otherwise zero; and $\Delta(i_0)\leq
0$  for $i_0> (tp+a)/q $ with exactly
$$B_t:=\#\{\gamma\not\in\Gamma: \gamma> (tp+a)/q\}$$ times taking the value $-1$,
otherwise zero. Notice that both $A_t$ and $B_t$ are strict
positive (since $0\in \Gamma$, respectively $2\delta-1\not\in
\Gamma$ and $2\delta-1>(tp+a)/q$). This shows that
\begin{equation*}
M_t:=\max_{0\leq i_0<m_f} \
\tau_a(tm_f+i_0)=\tau_a(tm_f)+A_t=\tau_a((t+1)m_f)+B_t,\tag{2}\end{equation*}
and
$$M_t>\max\{\, \tau_a(tm_f), \tau_a(tm_f+m_f)\,\}.$$
Therefore,  the graded root associated with the values
$\{\tau_a(i)\}_{0\leq i\leq (t_a+1)m_f}$ is the same as the graded
root associated with the values
$$\tau_a(0),M_0,\tau_a(m_f),M_1,\tau_a(2m_f), M_2, \ldots,
\tau_a(t_am_f), M_{t_a}, \tau_a(t_am_f+m_f).$$ Finally notice,
since $\#\{\gamma\not\in\Gamma\}=\delta$, one has
 $\delta-B_t=\#\{\gamma\not\in\Gamma: \gamma\leq (tp+a)/q\}$,
 hence
 $\delta-B_t+A_t=[(tp+a)/q]+1$.
 Hence, by (2),
 $$\tau_a((t+1)m_f)-\tau_a(tm_f)=\Big[\frac{tp+a}{q}\Big]+1-\delta.$$
Since $\tau_a(0)=0$, this gives $\tau_a(tm_f)$ inductively.
Notice also that $B_t=\alpha_{[(tp+a)/q]}$.

Clearly, the graded root associated with $\tau_a$ is the same as
the graded root associated with $\tilde{\tau}_a:\{0,1,2,\ldots,
2t_a+2\}\to \Z$, where $\tilde{\tau}_a(2t):=\tau_a(tm_f)$ and
$\tilde{\tau}_a(2t+1):=M_t$. This is the tau-function of
\ref{Maintheor} if we delete  \ $\tilde{}$\ .

\section{Examples}

\subsection{Example.}\label{epq1} Assume $p=q=1$.
 In this case $M$ is integral
homology sphere; $a=0$ and $t_0=2\delta-2=\mu-2$. In particular,
the rank of $\Ker U$ is $\mu$. Moreover, $r_0=\delta(\delta-1)$
and $\tau_0(2t)=t(t-2\delta+1)/2$. The reader is invited to draw
the graded root and verify that
$$HF^+(-M)=\calt_0^+\oplus \calt_0(\alpha_{\delta-1})\oplus\
\bigoplus_{i=1}^{\delta-1} \, \calt_{i(i+1)}(\,
\alpha_{i-1+\delta}\,)^{\oplus 2};$$
$$\Ker U=\bigoplus_{i=0}^{\delta-1}\ \Z_{(i^2+i)}^2.$$

\subsection{Example.}\label{ep1} More generally, assume only that
$p=1$. The Heegaard Floer homology of the unique
$spin^c$-structure is:
$$HF^+(-M)=\calt_0^+\oplus \calt_0(\alpha_{\delta-1})^{\oplus
q}\oplus \bigoplus_{i=1}^{(\delta-1)q}\,
\calt_{([l/q]+1)(\{l/q\}q+l)}(\, \alpha_{\delta-1+\lceil
l/q\rceil}\,)^{\oplus 2}.$$

\subsection{Remark.} Apparently in \ref{ep1},
$HF^+(-M)$ contains less information than the
polynomial $P$, the above formula involves only the coefficients
$\alpha_j$ for $j\geq \delta-1$. But this is not the case.
Indeed, since the Alexander polynomial  $\Delta(t)$ is symmetric,
one gets that
$$t^{\mu-2}P(1/t)-P(t)=\frac{\delta(1+t^{\mu-1})-(1+t+\cdots
+t^{\mu-1})}{t-1},$$ hence $\alpha_{\mu-2-i}-\alpha_i$ is a
universal number. In particular, from $\{\alpha _j\}_{j\geq
\delta-1}$ one can recover $P$.

This shows that from $HF^+(-S^3_{-1/q}(K_f))$ one can recover
both the integer $q$ and the isotopy type of $K_f\subset S^3$. It
looks that similar result is valid  for general  surgery
coefficients as well.

\subsection{Remark.} In the above example it is striking a
$\Z_2$-symmetry of the graded root and of $HF^+(-M)$. We will
explain this for $a=0$.

The point is that if the canonical characteristic element has only
integral coefficients, then all the theory associated with the
integral cycles (the lattice $L$ in section 5) has a duality.
This is happening for example if $p=1$; or if $q=1$ and
$2\delta-2$ is divisible by $p$. In our case, in these situations,
the function $\tau_0$ is stable with respect to the $\Z_2$ action
$i\mapsto 2t_0+2-i$, i.e. $\tau_0(i)=\tau_0(2t_0+2-i)$. This
induces a symmetry of the root and of the Heegaard Floer
homology. [In singularity theory a resolution graph whose
canonical characteristic element has only integral coefficients is
called `numerical Gorenstein graph'.]

But, for general $p/q$, the graphs are not `numerical Gorenstein'.
Even if they are, for  general $a$, the symmetry may fail.

\subsection{Example.} Assume that $K_f\subset S^3$ is the torus
knot $(4,5)$ (i.e. $g=1$ and $(p_1,q_1)=(4,5)$).  In this case
$\delta=6$, $\mu=12$,  $\Gamma$ is generated by 4 and 5, hence
$$P(t)=6+5t+4t^2+3t^3+3t^4+3t^5+2t^6+t^7+t^8+t^9+t^{10}.$$
Then for some $(p/q)$-surgeries the corresponding graded roots are
presented in the next Figure. [Here we did not draw all the
vertices, only those ones which are either local minimums or
supremums of pairs of local minimums; in fact, exactly these ones
are given by the function $\tau_a$. Also, the roots should be
continued upward with one vertex in $\chi^{-1}(k)$ for each
$k\geq 1$.]

Recall that when we compute the grading of $HF^+$, the value of
$\tau_a$ is doubled, and then shifted by $r_a$. The corresponding
values of $r_a$ in these five cases are: $$30, \ 71/4, \ 49/4, \
\frac{(p+10)^2-p}{4p},\ 60.$$ In particular, in the first and
fifth case (i.e. when $p=1$) one has $d(-M,\sigma_0)=0$, as we
expected: the Heegaard Floer homology are written in \ref{epq1},
resp. \ref{ep1}. In the second and third cases (i.e. when $p/q=2$
and $a=0,1$) the $\Z[U]$-modules are:
$$\Big(\ \ \calt_{-18}^+\oplus \calt_{-16}(2)^{\oplus 2}\oplus
\calt_{-10}(1)^{\oplus 2}\oplus \calt_{0}(1)^{\oplus 2}\ \
\Big)\,[71/4];$$
$$\Big(\ \ \calt_{-12}^+\oplus \calt_{-12}(3) \oplus
\calt_{-8}(1)^{\oplus 2}\oplus \calt_{0}(1)^{\oplus 2}\ \
\Big)\,[49/4].$$ In the forth case it is
$$\calt_{-10}^+\oplus \calt_0(1)\, \Big[\,\frac{(p+1)^2-p}{4p}\,\Big].$$

\begin{picture}(400,400)(20,20)

\dashline{1}(40,400)(430,400) \put(20,380){\makebox(0,0){0}}
\dashline{1}(40,390)(430,390) \put(20,380){\makebox(0,0){0}}
\dashline{1}(30,380)(430,380) \put(20,380){\makebox(0,0){0}}
\dashline{1}(40,370)(430,370) \put(20,380){\makebox(0,0){0}}
\dashline{1}(40,360)(430,360) \put(20,380){\makebox(0,0){0}}
\dashline{1}(40,350)(430,350) \put(20,380){\makebox(0,0){0}}
\dashline{1}(40,340)(430,340) \put(20,380){\makebox(0,0){0}}
\dashline{1}(30,330)(430,330) \put(20,330){\makebox(0,0){-5}}
\dashline{1}(40,320)(430,320) \put(20,380){\makebox(0,0){0}}
\dashline{1}(40,310)(430,310) \put(20,380){\makebox(0,0){0}}
\dashline{1}(40,300)(430,300) \put(20,380){\makebox(0,0){0}}
\dashline{1}(40,290)(430,290) \put(20,380){\makebox(0,0){0}}
\dashline{1}(30,280)(430,280) \put(20,280){\makebox(0,0){-10}}
\dashline{1}(40,270)(430,270) \put(20,380){\makebox(0,0){0}}
\dashline{1}(40,260)(430,260) \put(20,380){\makebox(0,0){0}}
\dashline{1}(40,250)(430,250) \put(20,380){\makebox(0,0){0}}
\dashline{1}(40,240)(430,240) \put(20,380){\makebox(0,0){0}}
\dashline{1}(30,230)(430,230) \put(20,230){\makebox(0,0){-15}}
\dashline{1}(40,220)(430,220) \put(20,380){\makebox(0,0){0}}
\dashline{1}(40,210)(430,210) \put(20,380){\makebox(0,0){0}}
\dashline{1}(40,200)(430,200) \put(20,380){\makebox(0,0){0}}
\dashline{1}(40,190)(430,190) \put(20,380){\makebox(0,0){0}}
\dashline{1}(30,180)(430,180) \put(20,180){\makebox(0,0){-20}}
\dashline{1}(40,170)(430,170) \put(20,380){\makebox(0,0){0}}
\dashline{1}(40,160)(430,160) \put(20,380){\makebox(0,0){0}}
\dashline{1}(40,150)(430,150) \put(20,380){\makebox(0,0){0}}
\dashline{1}(40,140)(430,140) \put(20,380){\makebox(0,0){0}}
\dashline{1}(30,130)(430,130) \put(20,130){\makebox(0,0){-25}}
\dashline{1}(40,120)(430,120) \put(20,380){\makebox(0,0){0}}
\dashline{1}(40,110)(430,110) \put(20,380){\makebox(0,0){0}}
\dashline{1}(40,100)(430,100) \put(20,380){\makebox(0,0){0}}
\dashline{1}(40,90)(430,90) \put(20,380){\makebox(0,0){0}}
\dashline{1}(30,80)(430,80) \put(20,80){\makebox(0,0){-30}}
\dashline{1}(40,70)(430,70) \put(20,380){\makebox(0,0){0}}

\put(80,380){\circle*{3}} \put(100,380){\circle*{3}}
\put(90,390){\circle*{3}} \put(90,340){\circle*{3}}
\put(80,330){\circle*{3}} \put(100,330){\circle*{3}}
\put(90,300){\circle*{3}} \put(80,290){\circle*{3}}
\put(100,290){\circle*{3}} \put(90,270){\circle*{3}}
\put(80,260){\circle*{3}} \put(100,260){\circle*{3}}
\put(90,260){\circle*{3}} \put(80,240){\circle*{3}}
\put(100,240){\circle*{3}} \put(80,230){\circle*{3}}
\put(100,230){\circle*{3}} \put(90,400){\line(0,-1){140}}
\put(80,380){\line(1,1){10}} \put(100,380){\line(-1,1){10}}
\put(80,330){\line(1,1){10}} \put(100,330){\line(-1,1){10}}
\put(80,290){\line(1,1){10}} \put(100,290){\line(-1,1){10}}
\put(80,260){\line(1,1){10}} \put(100,260){\line(-1,1){10}}
\put(80,240){\line(1,2){10}} \put(100,240){\line(-1,2){10}}
\put(80,230){\line(1,3){10}} \put(100,230){\line(-1,3){10}}

\put(160,380){\circle*{3}} \put(180,380){\circle*{3}}
\put(170,390){\circle*{3}} \put(170,340){\circle*{3}}
\put(160,330){\circle*{3}} \put(180,330){\circle*{3}}
\put(170,320){\circle*{3}} \put(160,300){\circle*{3}}
\put(180,300){\circle*{3}} \put(170,290){\circle*{3}}
\put(170,400){\line(0,-1){110}} \put(160,380){\line(1,1){10}}
\put(180,380){\line(-1,1){10}} \put(160,330){\line(1,1){10}}
\put(180,330){\line(-1,1){10}} \put(160,300){\line(1,2){10}}
\put(180,300){\line(-1,2){10}}

\put(240,380){\circle*{3}} \put(260,380){\circle*{3}}
\put(250,390){\circle*{3}} \put(250,350){\circle*{3}}
\put(240,340){\circle*{3}} \put(260,340){\circle*{3}}
\put(240,320){\circle*{3}} \put(260,320){\circle*{3}}
\put(250,400){\line(0,-1){50}} \put(240,380){\line(1,1){10}}
\put(260,380){\line(-1,1){10}} \put(240,340){\line(1,1){10}}
\put(260,340){\line(-1,1){10}} \put(240,320){\line(1,3){10}}
\put(260,320){\line(-1,3){10}}

\put(400,380){\circle*{3}} \put(420,380){\circle*{3}}
\put(410,390){\circle*{3}} \put(410,340){\circle*{3}}
\put(400,330){\circle*{3}} \put(420,330){\circle*{3}}
\put(410,290){\circle*{3}} \put(400,280){\circle*{3}}
\put(420,280){\circle*{3}} \put(410,250){\circle*{3}}
\put(400,240){\circle*{3}} \put(420,240){\circle*{3}}
\put(410,210){\circle*{3}} \put(400,200){\circle*{3}}
\put(420,200){\circle*{3}} \put(410,180){\circle*{3}}
\put(400,170){\circle*{3}}\put(420,170){\circle*{3}}
\put(410,150){\circle*{3}} \put(400,140){\circle*{3}}
\put(420,140){\circle*{3}} \put(410,130){\circle*{3}}
\put(400,120){\circle*{3}} \put(410,120){\circle*{3}}
\put(420,120){\circle*{3}} \put(410,110){\circle*{3}}
\put(400,100){\circle*{3}} \put(420,100){\circle*{3}}
\put(400,90){\circle*{3}} \put(420,90){\circle*{3}}
\put(400,80){\circle*{3}} \put(410,80){\circle*{3}}
\put(420,80){\circle*{3}}

\put(410,400){\line(0,-1){320}} \put(400,380){\line(1,1){10}}
\put(420,380){\line(-1,1){10}} \put(400,330){\line(1,1){10}}
\put(420,330){\line(-1,1){10}} \put(400,280){\line(1,1){10}}
\put(420,280){\line(-1,1){10}} \put(400,240){\line(1,1){10}}
\put(420,240){\line(-1,1){10}} \put(400,200){\line(1,1){10}}
\put(420,200){\line(-1,1){10}} \put(400,170){\line(1,1){10}}
\put(420,170){\line(-1,1){10}} \put(400,140){\line(1,1){10}}
\put(420,140){\line(-1,1){10}} \put(400,120){\line(1,1){10}}
\put(420,120){\line(-1,1){10}} \put(400,100){\line(1,2){10}}
\put(420,100){\line(-1,2){10}} \put(400,90){\line(1,2){10}}
\put(420,90){\line(-1,2){10}} \put(400,80){\line(1,3){10}}
\put(420,80){\line(-1,3){10}}

\put(320,380){\circle*{3}} \put(330,390){\circle*{3}}
\put(330,330){\circle*{3}} \put(320,380){\line(1,1){10}}
\put(330,400){\line(0,-1){70}}

\put(330,50){\makebox(0,0)[t]{$p\geq 11,\ q=1 $}}
\put(330,35){\makebox(0,0)[t]{$a=0$}}
\put(90,50){\makebox(0,0)[t]{$p/q=1$}}
\put(90,35){\makebox(0,0)[t]{$a=0$}}
\put(170,50){\makebox(0,0)[t]{$p/q=2$}}
\put(170,35){\makebox(0,0)[t]{$a=0$}}
\put(250,50){\makebox(0,0)[t]{$p/q=2$}}
\put(250,35){\makebox(0,0)[t]{$a=1$}}
\put(410,50){\makebox(0,0)[t]{$p/q=1/2$}}
\put(410,35){\makebox(0,0)[t]{$a=0$}}
\put(20,400){\makebox(0,0)[t]{$\tau_a(i)$}}

\end{picture}

\end{document}